\documentclass[12pt, reqno]{amsart}
\usepackage{amssymb,amsthm,amsmath}
\usepackage[numbers,sort&compress]{natbib}
\usepackage{amssymb,amsmath}
\usepackage{amsfonts}
\usepackage{mathrsfs}
\usepackage{latexsym}
\usepackage{amssymb}
\usepackage{amsthm}
\usepackage{pdfsync}
\usepackage{indentfirst}
\date{June 18, 2021}

\usepackage{amsmath}
\usepackage{amsthm}

\newtheorem{theorem}{Theorem}

\newtheorem{remark}[theorem]{Remark}

\newtheorem{lemma}[theorem]{Lemma}

\newcommand{\beq}{\begin{equation}}
\newcommand{\eeq}{\end{equation}}
\newcommand{\ben}{\begin{eqnarray}}
\newcommand{\een}{\end{eqnarray}}
\newcommand{\beno}{\begin{eqnarray*}}
\newcommand{\eeno}{\end{eqnarray*}}

\numberwithin{equation}{section}

\allowdisplaybreaks
\begin{document}

\title[Liouville type theorems for the steady MHD and Hall-MHD equations]{Remarks on Liouville type theorems for the steady MHD and Hall-MHD equations}

\author{Xiaomeng~Chen}
\address[Xiaomeng~Chen]{School of Mathematical Sciences, Dalian University of Technology, Dalian, 116024,  China}
\email{cxm2381033@163.com}

\author{Shuai~Li}
\address[Shuai~Li]{School of Mathematical Sciences, Dalian University of Technology, Dalian, 116024,  China}
\email{ls21701048@163.com}

\author{Wendong~Wang}
\address[Wendong~Wang]{School of Mathematical Sciences, Dalian University of Technology, Dalian, 116024,  China}
\email{wendong@dlut.edu.cn}




\begin{abstract}
In this note we investigate Liouville type theorems for the steady  three dimensional MHD and Hall-MHD equations, and show that the velocity field $u$ and the magnetic field $B$ are vanishing provided that $B\in L^{6,\infty}(\mathbb{R}^3)$ and $u\in BMO^{-1}(\mathbb{R}^3)$, which state that the velocity field plays an important role. Moreover, the similar result holds in the case of partial viscosity or diffusivity for the three dimensional MHD equations.

\end{abstract}

\maketitle

\allowdisplaybreaks

\section{Introduction}

Consider the following general system including the so called magnetohydrodynamics (MHD) and Hall-MHD equations:
\begin{eqnarray}\label{eq:MHD}
 \left\{
    \begin{array}{llll}
    \displaystyle - \Delta u + u \cdot \nabla u + \nabla p = B \cdot \nabla B, ~~ \rm{in} ~~ \mathbb{R}^3, \\
    \displaystyle - \Delta B + u \cdot \nabla B - B \cdot \nabla u = \alpha \nabla \times ((\nabla \times B) \times B), ~~ \rm{in} ~~ \mathbb{R}^3,\\
    \displaystyle {\rm{div}} ~ u = 0, ~~~ {\rm{div}} ~ B = 0, ~~ \rm{in} ~~ \mathbb{R}^3,
    \end{array}
 \right.
\end{eqnarray}
where $u = u(x)=(u_1,u_2,u_3)$ is the velocity field of the fluid flows, $B = B(x) = (B_1,B_2,B_3)$ is the magnetic field,  $p = p(x)$ is the pressure of the flows, and $\alpha \in \mathbb{R}$ is an absolute constant. When $\alpha = 0$,  the system (\ref{eq:MHD}) is the MHD system, which describes the steady state of the magnetic properties of electrically conducting fluids, including
 plasmas, liquid metals, etc;  for the physical background
and mathematical theory we refer to Schnack \cite{Sch} and the references therein. When $\alpha = 1$,  this system governs the dynamics plasma flows of strong shear magnetic fields as in the solar flares, and has many important applications in the astrophysics (for example, see Chae-Degond-Liu \cite{CDL}).

In this note, we focus on the Liouville type properties of the above system, which is motivated by the development of Navier-Stokes equations.
When $B = 0$,  the system (\ref{eq:MHD}) reduces to the standard Navier-Stokes equations, and
a very
challenging open question is if there exists a nontrivial solution when the Dirichlet integral $\int_{\mathbb{R}^3}|\nabla u|^2dx$ is finite, which
dates back to Leray's celebrated paper \cite{Leray} and is explicitly written in Galdi's book (\cite{Galdi}, Remark X. 9.4, pp.729; see also Tsai's book \cite{Tsai-2018}, pp.23).
This uniqueness problem,
or equivalently the Liouville type problem is widely open. Galdi proved the above Liouville type theorem by assuming $u\in L^{\frac92}(\mathbb{R}^3)$ in \cite{Galdi}.
Chae in \cite{Chae2014} showed the condition $\triangle u\in L^{\frac65}(\mathbb{R}^3)$ is sufficient for the vanishing property of $u$ by exploring the maximum principle of the head pressure. Also, Chae-Wolf gave an improvement of logarithmic form for Galdi's result in \cite{ChaeWolf}
by assuming that $\int_{\mathbb{R}^3} |u|^{\frac92}\{\ln(2+\frac{1}{|u|})\}^{-1}dx<\infty$.
 Seregin obtained the conditional criterion $u\in BMO^{-1}(\mathbb{R}^3)$ in \cite{Se}.
More references, we refer to \cite{KPR,KTW,Se2,SeWang,CPZZ,Chae-Wolf2019,Chae2021} and the references therein.
Relatively speaking, the two-dimensional case is more easier due to the vorticity of the 2D  Navier-Stokes equations satisfies a nice elliptic equation,
to which the maximum principle applies. For example, Gilbarg-Weinberger \cite{GW1978} obtained Liouville type theorem provided the Dirichlet energy is finite. As a different type of Liouville property for the 2D Navier-Stokes equations,
Koch-Nadirashvili-Seregin-Sverak \cite{KNSS} showed
that any bounded solution is trivial solution, say $u \equiv C$,
as a byproduct of their results on the non-stationary case (see also \cite{BFZ2013} for the unbounded velocity).

However, for the MHD tpye system like (\ref{eq:MHD}), the situation is quite different. Due to the lack of maximum principle, there is not much progress in the study of MHD equation. For the two-dimensional MHD equations,  Liouville type theorems were proved by assuming the smallness of the norm of the magnetic field in \cite{WW} and \cite{Wang}, and we refer to the recent paper in \cite{DHS} by removing the smallness assumption. For the three-dimensional MHD equations,
Chae-Weng \cite{ChaeWeng} proved that if a smooth solution to (\ref{eq:MHD}) in $\mathbb{R}^3$ with finite Dirichlet integral
\beno\label{eq: energy}
\int_{\mathbb{R}^3} |\nabla u|^2 + |\nabla B|^2 dx < +\infty,
\eeno
and $u \in L^3(\mathbb{R}^3)$, then the solutions $(u,B)$ is identically zero.
Schulz \cite{Schulz} proved that if a smooth solution of the stationary MHD equations is in $L^6(\mathbb{R}^3)$ and $u,B \in BMO^{-1}(\mathbb{R}^3)$, then it is identically zero. Recently, Chae-Wolf \cite{Chae-Wolf-MHD} show that $L^6$ mean oscillations of the potential function of the velocity and magnetic field have certain linear growth by using the technique of \cite{Chae-Wolf2019}. More references, we refer to \cite{ZYQ,LLN,LP,YX,FW} and the references therein.

Motivated by some numerical experiments in \cite{PP}, which seem to indicate that the velocity field should play a more important role than the magnetic field in the regularity theory of solutions to the MHD equations. Our goal is to get rid of the extra magnetic conditions in \cite{ChaeWeng} or \cite{Schulz}.
Recall the definition of the space $BMO^{-1}(\mathbb{R}^3)$ as follows. For a measurable set $E \subset \mathbb{R}^3$, we denote by $|E|$ the 3-dimensional Lebesgue measure of $E$, and for $f \in L^1(E)$ we use the notation
\ben\label{eq:f_mean }
f_E = \frac{1}{|E|} \int_E f dx,
\een
Combining (\ref{eq:f_mean }) with H\"{o}lder's inequality, we know for all $f\in L^p(E)$
\beno\label{ine:age l p}
||f - f_E||_{L^p(E)} \leq C ||f||_{L^p(E)}.
\eeno
Then we say that $u \in BMO^{-1}(\mathbb{R}^3)$, if there exists $\mathbf{\Phi} \in C^\infty(\mathbb{R}^3,\mathbb{R}^{3 \times 3})$ such that $\nabla \cdot \mathbf{\Phi} = u$, and for all $r>0$, there holds
\beno
\sup_{r>0}\frac{1}{|B(r)|} \int_{B(r)} |\mathbf{\Phi} - \mathbf{\Phi}_{B(r)}| dx <\infty,
\eeno
i.e. $\mathbf{\Phi}\in  BMO(\mathbb{R}^3)$.
And a remarkable property of BMO space is
\ben\label{eq: BMO property}
\sup_{r>0}\frac{1}{|B(r)|} \int_{B(r)} |\mathbf{\Phi} - \mathbf{\Phi}_{B(r)}|^p dx <\infty, \quad \forall ~~ 1 \leq  p < +\infty.
\een
Next, the definition of the space $L^{p,q}(\mathbb{R}^n)$ is as follow:
Let $\Omega \subset \mathbb{R}^n$ and $1\leq p,l\leq \infty$, we say that a measurable function $f$ is in $L^{p,l}(\Omega)$ if $\|f\|_{L^{p,l}(\Omega)}<+\infty$, where
\begin{equation}\nonumber
\|f\|_{L^{p,l}}(\Omega)=
\left\{\begin{array}{llll}
\left(\int^\infty_0 \lambda^{l-1}|\{x\in \Omega;|f|>\lambda\}|^\frac{l}{p} d\lambda \right)^\frac1l,\quad ~{\rm for}~r<\infty;\\
sup_{\lambda>0} \lambda|\{x\in\Omega;|f|>\lambda\}|^\frac1p, \quad ~{\rm for}~r=\infty.
\end{array}\right.
\end{equation}
We use H\"{o}lder inequality in Lorentz spaceas follows (\cite{Neil}).
Assume $1 \leq p_1,p_2 \leq \infty$, $1 \leq q_1,q_2 \leq \infty$ and $u \in L^{p_1,q_1}(\Omega)$, $v \in L^{p_2,q_2}(\Omega)$. Then $uv \in L^{p,q}(\Omega)$ and $\frac 1q \leq \frac{1}{q_1} + \frac{1}{q_2}$ and the inequality
\ben\label{ine:Holder}
||uv||_{L^{p,q}(\Omega)} \leq C ||u||_{L^{p_1,q_1}(\Omega)}||v||_{L^{p_2,q_2}(\Omega)}
\een
is valid,where $\frac 1p = \frac{1}{p_1} + \frac{1}{p_2}$.

Our main result is as follows:

\begin{theorem}\label{thm:lorentz}
Let $(u,B,p) \in C^\infty(\mathbb{R}^3) \times C^\infty(\mathbb{R}^3) \times C^\infty(\mathbb{R}^3)$ be a solution of (\ref{eq:MHD}) with $B \in L^{6,\infty}(\mathbb{R}^3)$. Assume that $u \in BMO^{-1}(\mathbb{R}^3)$.
Then $u = B = 0$.
\end{theorem}

\begin{remark}\label{rem:lorentz} The above theorem shows that the velocity field should play a more important role than the magnetic field in the uniqueness theory, which is similar as in the regularity theory.
For example, \cite{HX1} and  \cite{WZ}  have presented some regularity criterions to
the MHD equations in terms of the velocity field only.
\end{remark}

\begin{remark}\label{rem:lorentz2}
The above result generalized Seregin's result  in \cite{Se} or \cite{Se2} to the MHD case or Hall-MHD case when $B\equiv0$, which also improved Chae-Weng's result \cite{ChaeWeng}, where they assumed $\nabla u,\nabla B\in L^2(\mathbb{R}^3)$ and $u\in L^3(\mathbb{R}^3)$. It also improved Schulz's theorem \cite{Schulz} by removing the additional conditions $B\in BMO^{-1}(\mathbb{R}^3)$ and $u \in L^6(\mathbb{R}^3)$. Moreover, we relaxed the condition of $B$ to the Lorentz space of $L^{6,\infty}(\mathbb{R}^3)$. A new observation is the use of delicate $L^q$ estimates of stationary Stokes system, which is independent of interest.
\end{remark}

For the MHD equations with partial viscosity or diffusivity, we have the following similar conclusions.
Consider the MHD system (\ref{eq:MHD}) as follows:
\begin{eqnarray}\label{eq:VMHD}
 \left\{
    \begin{array}{llll}
    \displaystyle - \lambda_1 \partial_{11} u - \lambda_2 \partial_{22} u - \lambda_3 \partial_{33} u + u \cdot \nabla u + \nabla p = B \cdot \nabla B, ~~ \rm{in} ~~ \mathbb{R}^3, \\
    \displaystyle - \mu_1 \partial_{11} B - \mu_2 \partial_{22} B - \mu_3 \partial_{33} B + u \cdot \nabla B = B \cdot \nabla u, \\
    \displaystyle {\rm{div}} ~ u = 0, ~~~ {\rm{div}} ~ B = 0, \\
    \displaystyle \lambda_i \geq 0, ~~~ \mu_i \geq 0, i = 1,2,3. \\
    \end{array}
 \right.
\end{eqnarray}

For the system (\ref{eq:VMHD}), we have
\begin{theorem}\label{thm:viscidity}
Let $u,b$ be a smooth solutions of the system (\ref{eq:VMHD}). Let's further assume that $B \in L^6(\mathbb{R}^3)$ and $u \in L^3(\mathbb{R}^3)$, then $u, B \equiv 0$ if the following condition holds
\ben\label{eq:condition}
\lambda_1 + \lambda_2 + \lambda_3 > 0, ~~~ \mu_1 + \mu_2 + \mu_3 > 0.
\een
\end{theorem}

\begin{remark}\label{rem:lorentz3}
When $\lambda_1 + \lambda_2 + \lambda_3=0$, but $\mu_1 + \mu_2 + \mu_3 > 0$, one can derive that $B\equiv0$ and $u$ satisfies the three dimensional Euler equations.  It is impossible to deduce that $u$ is vanishing from the  bounded-ness of $L^q$ norm. For example, we can refer to the counterexample belonging to $C_0^\infty(\mathbb{R}^3)$ in \cite{Ga}. The same is true for another situation of $\mu_1 + \mu_2 + \mu_3 = 0$.
\end{remark}


%
%

In the proof, we need the following lemma (see, for example, Lemma A.5 \cite{WZ1}).
%
%

\begin{lemma}\label{lem:WZ}
Let $f, g \in \dot{W}^{1,2}(\mathbb{R}^3) \cap BMO^{-1}(\mathbb{R}^3)$. Then there holds
\beno
||fg||_{L^2} \leq C \left(||f||_{\dot{W}^{1,2}}||g||_{BMO^{-1}} + ||f||_{BMO^{-1}}||g||_{\dot{W}^{1,2}}\right).
\eeno
\end{lemma}

We also need an interpolation inequality in Lorentz space (see, for example, Theorem 2.1 \cite{NJN}).
\begin{lemma}\label{lem:NJN}
Let $\Omega$ be a domain in $\mathbb{R}^n$. Let $0 < q < p < r \leq +\infty$ and $\alpha > 0$. If $f \in L^{q,\infty}(\Omega) \cap L^{r,\infty}(B(r))$, then $f \in L^{p,\alpha}(\Omega)$ and
\beno
||f||_{L^{p,\alpha}(\Omega)} \leq C ||f||_{L^{q,\infty}(\Omega)}^\theta ||f||_{L^{r,\infty}(\Omega)}^{1-\theta},
\eeno
where $C = C(q,r,p,\alpha) > 0$ and
\beno
\frac 1p = \frac{\theta}{q} + \frac{1-\theta}{r}.
\eeno
\end{lemma}

\begin{remark}
The lemma \ref{lem:NJN} follows from Theorem 2.1 in \cite{NJN}, where the domain is the whole space $\mathbb{R}^n$. Note that $||f\chi_\Omega||_{L^{q,\infty}(\mathbb{R}^n)} = ||f||_{L^{q,\infty}(\Omega)}$, which show that the lemma holds for the general domain $\Omega$  in $\mathbb{R}^n$.
\end{remark}

Throughout this article, $C(\alpha_1,\cdots,\alpha_n)$ denotes a constant depending on $\alpha_1,\cdots,\alpha_n$, which may be different from line to line.

\section{Proof of Theorem \ref{thm:lorentz}}

Before starting the proof of Theorem \ref{thm:lorentz}, we state the following interior $L^q$ estimates for the steady Stokes system of an alternative version and sketch its proof.
\begin{lemma}\label{lem:Tsai's}
Let $B(r) \subset B(R)\subset \mathbb{R}^n$ with $n\geq 2$ be concentric balls with $0 < r < R$. Assume that $v$ is a very weak solution of the following Stokes system
\ben\label{eq:stokes}
- \Delta v_i + \partial_i p = g_i, \quad {\rm div}\,v = 0 \quad {\rm in}~ B(R),
\een
where $g_i \in L^q(B(R))$, $1 < q < \infty$. Then  $v \in W_{loc}^{2,q}(B(R))$ and there exists  a function $p\in W_{loc}^{1,q}(B(R))$. Moreover, there holds
\ben\label{eq:stokes estimate}
||\nabla^2 v||_{L^q(B(r))} + ||\nabla p||_{L^q(B(r))}\leq C(n,q) \left(||g||_{L^q(B(R))} +  (R-r)^{-2}||v||_{L^q(B(R))}\right).
\een
\end{lemma}

\begin{remark}
A very weak solution $v$ of the system of (\ref{eq:stokes}) in $\Omega\subset\mathbb{R}^n$ is defined as follows. The divergence-free vector field $v\in L_{loc}^{1}(\Omega)$ satisfies
\beno
\int v\cdot \triangle\zeta=-\int g\cdot \zeta,\quad \forall~\zeta\in C_{c,\sigma}^\infty(\Omega),
\eeno
where $ C_{c,\sigma}^\infty(\Omega)$ denotes any times differentiable divergence-free vector fields with compact
support in $\Omega.$
\end{remark}

\begin{proof}
At first, recall
Lemma 2.12 in \cite{Tsai-2018} and we have
\beno\label{eq:stokes estimate-Tsai}
||\nabla^2 v||_{L^q(B(r))} \leq C(n,q,r,R) \left(||g||_{L^q(B(R))} +  ||v||_{L^1(B(R))}\right).
\eeno
Specially,
\beno
||\nabla^2 v||_{L^q(B(1))} \leq C(n,q) \left(||g||_{L^q(B(2))} +  ||v||_{L^1(B(2))}\right).
\eeno
Hence, for any $0<r<R$ and $2r\leq R$, by scaling we get
\ben\label{eq:stokes estimate-Tsai2}
||\nabla^2 v||_{L^q(B(r))} \leq C(n,q) \left(||g||_{L^q(B(2r))} + r^{-n-2+\frac{n}{q}} ||v||_{L^1(B(2r))}\right).
\een
Secondly, for any fixed $r>0$ with $r<R$, and $x_0 \in B(r)$, let $10d = R - r>0$, then
\beno
\bigcup_{x_0\in B(r)}B(x_0,d)\supset B(r),
\eeno
where $B(x_0,d)$ denotes the ball of radius $d$ centered at $x_0.$ Due to the lemma of Vitali, there exists finite points, $x_1,\cdots,x_k\in B(r)$ such that $B(x_1,d),\cdots, B(x_k,d)$  are  disjoint balls satisfying
\beno
B(r)\subset \bigcup_{j=1}^k B(x_j,5d),
\eeno
and at most $K_0$ finite balls with $K_0\leq 6^n$ have a joint point in $ B(R)$.
Thence, by (\ref{eq:stokes estimate-Tsai2}) and H\"{o}lder inequality we have
\beno
||\nabla^2 v||_{L^q(B(r))}^q
 &\leq& \sum_{j=1}^k ||\nabla^2 v||_{L^q(B(x_j,5d))}^q  \\
&\leq& C(n,q) \sum_{j=1}^k \left(||g||_{L^q(B(x_j,10d))}^q + d^{-(n+2)q+n} ||v||_{L^1(B(x_j,10d))}^q\right)\\
&\leq& C(n,q) \sum_{j=1}^k \left(||g||_{L^q(B(x_j,10d))}^q + d^{-2q} ||v||_{L^q(B(x_j,10d))}^q\right)
\eeno
Note that
at most $K_1$ finite balls of $B(x_j,10d)$ with $K_1\leq 11^n$ have a joint point in $ B(R)$. Then
\beno
||\nabla^2 v||_{L^q(B(r))}^q &\leq& C(n,q)11^n\left(||g||_{L^q(B(R))}^q + d^{-2q} ||v||_{L^q(B(R))}^q\right)\\
&\leq& C(n,q)\left(||g||_{L^q(B(R))}^q + (R-r)^{-2q} ||v||_{L^q(B(R))}^q\right)
\eeno
by using $10d = R - r$ and
\beno
B(R)\supset \bigcup_{j=1}^k B(x_j,10d).
\eeno

At last, using the equations of (\ref{eq:stokes}), we get
\beno
||\nabla p ||_{L^q(B(r))}
&\leq& C(n,q)\left(||g||_{L^q(B(R))}+ (R-r)^{-2} ||v||_{L^q(B(R))}\right)
\eeno
The proof is complete.
\end{proof}

Next we prove Theorem \ref{thm:lorentz}.
\begin{proof}
{\bf Step I. Cacciopolli type inequality.}

First, we introduce
some cut-off functions of $\phi_i$ with $i=1,\cdots,3$ for $r>1$ and $\rho, \tau>0$, which will be used later.
Assume that $\phi_i(x) \in C_0^{\infty}(\mathbb{R}^3)$ with $ 0 \le \phi_i \le 1$ for $i=1,\cdots,3$ satisfying the following conditions:\\
i) $r \le \rho < \tau \le 2r$;\\
ii)
\begin{equation}\label{eq:phi 1}
\phi_1(x)=\phi_1(|x|)=
\left\{\begin{array}{llll}
1,\quad ~{\rm in}~B_{\rho};\\
0, \quad ~{\rm in} ~ \mathbb{R}^3 \backslash B(\frac{\tau + 2\rho}{3});\\
\end{array}\right.
\end{equation}
iii)\begin{equation}\label{eq:psi 2}
\phi_2(x)=\phi_2(|x|)=
\left\{\begin{array}{llll}
1,\quad ~{\rm in}~B(\frac{\tau + 2\rho}{3});\\
0, \quad ~{\rm in} ~\mathbb{R}^3 \backslash B(\frac{2\tau + \rho}{3});\\
\end{array}\right.
\end{equation}
iv)\begin{equation}\label{eq:psi 3}
\phi_3(x)=\phi_3(|x|)=
\left\{\begin{array}{llll}
1,\quad ~{\rm in}~B(\frac{2\tau + \rho}{3});\\
0, \quad ~{\rm in} ~\mathbb{R}^3 \backslash B(\tau);\\
\end{array}\right.
\end{equation}
v)
\beno
|\nabla \phi_i| \leq \frac{C}{\tau - \rho},\quad |\nabla^2 \phi_i| \leq \frac{C}{(\tau - \rho)^2},\quad  i=1,\cdots,3.
\eeno

Secondly,
multiply $(\ref{eq:MHD})_1$ by $u \phi_1$, $(\ref{eq:MHD})_2$ by $(B - b) \phi_1$, and integrate them over $B(2r)$. Then by integration by parts we obtain
\beno
&&\int_{B(2r)} |\nabla u|^2 \phi_1 + |\nabla B|^2 \phi_1 dx \\
&=& \frac 12 \int_{B(2r)} |u|^2 \Delta \phi_1 dx + \frac12\int_{B(2r)} |u|^2 u \cdot \nabla \phi_1 dx \\
&&+ \int_{B(2r)} (p - a) u \cdot \nabla \phi_1 dx + \frac 12 \int_{B(2r)} |B - b|^2 \Delta \phi_1 dx \\
&&+ \frac12\int_{B(2r)} |B - b|^2 u \cdot \nabla \phi_1 dx - \int_{B(2r)} u \cdot (B - b) B \cdot \nabla \phi_1 dx \\
&&- \alpha \int_{B(2r)} ((\nabla \times (B - b)) \times B) \cdot ((B - b) \times \nabla \phi_1) dx\\
&\doteq& I_1 + \cdots+ I_7,
\eeno
where  $a$ and $b$ are constant vectors in $B(2r)$, to be decided, and we used the following formula for the last term of the right hand
\beno
\int \nabla \times F \cdot G\phi dx&=&\int \nabla \times G \cdot F\phi dx+\int G \times F \cdot \nabla \phi dx\\
&=&\int \nabla \times G \cdot F \phi dx - \int F \cdot (G \times  \nabla \phi) dx.
\eeno

For the terms of $I_1$ and $I_2$, it's easy to deduce that
\ben\label{ine:u 2}
|I_1| \leq C (\tau - \rho)^{-2} \int_{B(\frac{\tau + 2\rho}{3}) \setminus B(\rho)} |u|^2 dx.
\een
and
\ben\label{ine:u 3}
|I_2| \leq C (\tau - \rho)^{-1} \int_{B(\frac{\tau + 2\rho}{3}) \setminus B(\rho)} |u|^3 dx.
\een

For the term of $I_4$, by (\ref{ine:Holder}) we have
\ben\label{eq:I4}
|I_4| &\leq& C (\tau - \rho)^{-2} \int_{B(\frac{\tau + 2\rho}{3}) \setminus B(\rho)} |B - b|^2 dx\nonumber\\
&\leq & C (\tau - \rho)^{-2} ||B - b||_{L^{6,\infty}(B(\frac{2\tau + \rho}{3})\backslash B(\rho))}^2|| 1 ||_{L^{\frac32,1}(B(2r))}\nonumber\\
&\leq & C r^2 (\tau - \rho)^{-2} ||B - b||_{L^{6,\infty}(B(\frac{2\tau + \rho}{3})\backslash B(\rho))}^2.
\een

For the term $I_5$, using $\nabla \cdot \mathbf{\Phi} = u$, integration by parts and H\"{o}lder inequality, we have
\beno
|I_5| &=& \left|\int_{B(\frac{\tau + 2\rho}{3})\backslash B(\rho)} |B - b|^2 \partial_i(\Phi_{ij} - \bar{\Phi}_{ij}) \partial_j \phi_1 dx\right| \\
&\leq& C (\tau - \rho)^{-1} \int_{B(\frac{\tau + 2\rho}{3}) \backslash B(\rho)} |B - b| |\nabla B| |\mathbf{\Phi} - \bar{\mathbf{\Phi}}| dx\\
 &&+ C (\tau - \rho)^{-2} \int_{B(\frac{\tau + 2\rho}{3}) \backslash B(\rho)} |B - b|^2 |\mathbf{\Phi} - \bar{\mathbf{\Phi}}| dx \\
&\leq& C (\tau - \rho)^{-1} ||B-b||_{L^{6,\infty}(B(\frac{\tau + 2\rho}{3}) \backslash B(\rho))} ||\nabla B||_{L^2(B(\frac{\tau + 2\rho}{3}) \backslash B(\rho))} ||\mathbf{\Phi} - \bar{\mathbf{\Phi}}||_{L^{3,2}(B(\frac{\tau + 2\rho}{3}) )} \\
&&+ C (\tau - \rho)^{-2} ||B - b||_{L^{6,\infty}(B(\frac{\tau + 2\rho}{3}) \backslash B(\rho))}^2 ||\mathbf{\Phi} - \bar{\mathbf{\Phi}}||_{L^{\frac32,1}(B(\frac{\tau + 2\rho}{3}))}.
\eeno
Choosing $\bar{\mathbf{\Phi}} = \mathbf{\Phi}_{B(2r)}$ and  $b=B_{B(\frac{4r}{3})\backslash B(r)}$, by Lemma \ref{lem:NJN} and (\ref{eq: BMO property}), we have
\ben\label{eq:NJN1} \nonumber
||\mathbf{\Phi} - \bar{\mathbf{\Phi}}||_{L^{3,2}(B(2r))}
&\leq& C ||\mathbf{\Phi} - \bar{\mathbf{\Phi}}||_{L^2(B(2r))}^\frac 13  ||\mathbf{\Phi} - \bar{\mathbf{\Phi}}||_{L^4(B(2r))}^\frac 23 \leq C r,
\een
\ben\label{eq:NJN2} \nonumber
||\mathbf{\Phi} - \bar{\mathbf{\Phi}}||_{L^{\frac 32,1}(B(2r))}
&\leq& C ||\mathbf{\Phi} - \bar{\mathbf{\Phi}}||_{L^\frac 54(B(2r))}^\frac{5}{12}  ||\mathbf{\Phi} - \bar{\mathbf{\Phi}}||_{L^\frac 74(B(2r))}^\frac{7}{12} \leq C r^2,
\een
and by $B \in L^{6,\infty}(\mathbb{R}^3)$
\ben\label{eq:estimate b}
||b||_{L^{6,\infty}(B(\frac{\tau + 2\rho}{3}) \backslash B(\rho)}&\leq &||b||_{L^{6,\infty}(B(2r))}\leq Cr^{\frac12}|b|\nonumber\\
&\leq &Cr^{\frac12-3}||B||_{L^{6,\infty}(B(2r))}||1||_{L^{\frac65,1}(B(2r))}\leq C.
\een
Then we get
\ben\label{eq:I5}
|I_5| &\leq& C r (\tau - \rho)^{-1} ||\nabla B||_{L^2(B(\frac{2\tau + \rho}{3})\backslash B(\rho))} \nonumber\\
&&+ C r^2 (\tau - \rho)^{-2} ||B - b||_{L^{6,\infty}(B(\frac{2\tau + \rho}{3})\backslash B(\rho))}^2.
\een

The term of $I_6$ can be controlled by
\beno
|I_6| &=& \left|\int_{B(\frac{\tau + 2\rho}{3})\backslash B(\rho)} \partial_i(\Phi_{ij} - \bar{\Phi}_{ij}) (B - b)_j B_k \partial_k \phi_1 dx\right| \\
&\leq& C (\tau - \rho)^{-1} \int_{B(\frac{\tau + 2\rho}{3})\backslash B(\rho)} |\mathbf{\Phi} - \bar{\mathbf{\Phi}}| |B - b| |\nabla B| dx \\
 &&+ C (\tau - \rho)^{-1} \int_{B(\frac{\tau + 2\rho}{3})\backslash B(\rho)} |\mathbf{\Phi} - \bar{\mathbf{\Phi}}| |B| |\nabla B| dx \\
 &&+ C (\tau - \rho)^{-2} \int_{B(\frac{\tau + 2\rho}{3})\backslash B(\rho)} |\mathbf{\Phi} - \bar{\mathbf{\Phi}}| |B - b| |B| dx.
\eeno
Similarly as $I_5$, we have
\ben\label{eq:I6}
|I_6| &\leq&C r (\tau - \rho)^{-1} ||\nabla B||_{L^2(B(\frac{\tau + 2\rho}{3})\backslash B(\rho))} \nonumber\\
&& + C r^2 (\tau - \rho)^{-2} ||B - b||_{L^{6,\infty}(B(\frac{\tau + 2\rho}{3})\backslash B(\rho))}.
\een

Using (\ref{ine:Holder}) and (\ref{eq:estimate b}), the term of $I_7$ can be controlled by
\ben\label{eq:I7}
|I_7|
 &\leq&  C (\tau - \rho)^{-1} \int_{B(\frac{\tau + 2\rho}{3})\backslash B(\rho)} |\nabla B| |B| |B - b| dx \nonumber\\
&\leq& C (\tau - \rho)^{-1} ||\nabla B||_{L^2(B(\frac{\tau + 2\rho}{3})\backslash B(\rho))} ||B||_{L^{6,\infty}(B(\frac{\tau + 2\rho}{3}))} ||B - b||_{L^{6,\infty}(B(\frac{\tau + 2\rho}{3}))} ||1||_{L^{6,2}(B(\frac{\tau + 2\rho}{3}))}\nonumber \\
&\leq&  C r^\frac 12 (\tau - \rho)^{-1} ||\nabla B||_{L^2(B(\frac{\tau + 2\rho}{3})\backslash B(\rho))}.
\een
At last, we estimate the term of $I_3$, and by $\nabla \cdot \mathbf{\Phi} = u$ we get
\beno
|I_3| &=& \left|\int_{B(\frac{\tau + 2\rho}{3})\backslash B(\rho)} (p - a) \partial_i(\Phi_{ij} - \bar{\Phi}_{ij}) \partial_j \phi_1 dx\right| \\
&=& \left|\int_{B(\frac{\tau + 2\rho}{3})\backslash B(\rho)} (p - a) (\Phi_{ij} - \bar{\Phi}_{ij}) \partial_{ij} \phi_1 dx + \int_{B(\frac{\tau + 2\rho}{3})\backslash B(\rho)} \partial_ip (\Phi_{ij} - \bar{\Phi}_{ij}) \partial_j \phi_1 dx \right| \\
&\leq& C (\tau - \rho)^{-2} ||p - a||_{L^{\frac{3s}{3-s}}(B(\frac{\tau + 2\rho}{3})\backslash B(\rho))} ||\mathbf{\Phi} - \mathbf{\bar{\Phi}}||_{L^{\frac{3s}{4s-3}}(B(\frac{\tau + 2\rho}{3})\backslash B(\rho))} \\
&&+ C (\tau - \rho)^{-1} ||\nabla p||_{L^s(B(\frac{\tau + 2\rho}{3})\backslash B(\rho))} ||\mathbf{\Phi} - \mathbf{\bar{\Phi}}||_{L^{\frac{s}{s-1}}(B(\frac{\tau + 2\rho}{3})\backslash B(\rho))} \\
&\leq& C r^{\frac{4s-3}{s}}(\tau - \rho)^{-2} ||p - a||_{L^{\frac{3s}{3-s}}(B(\frac{\tau + 2\rho}{3})\backslash B(\rho))} +C r^{\frac{3s-3}{s}}(\tau - \rho)^{-1} ||\nabla p||_{L^s(B(\frac{\tau + 2\rho}{3})\backslash B(\rho))},
\eeno
where $s \in (1, 3)$.

At last, collecting (\ref{ine:u 2})-(\ref{eq:I4}), (\ref{eq:I5})-(\ref{eq:I7}) and the estimates of $I_3$, noting $r>1$ and $||B - b||_{L^{6,\infty}(B(2r))}\leq C$ from (\ref{eq:estimate b})  we arrive at the following Cacciopolli type inequality
\ben\label{ine:cacciopolli} \nonumber
&&\int_{B(\rho)} |\nabla u|^2 + |\nabla B|^2 dx \\ \nonumber
&&\leq C (\tau - \rho)^{-2} \int_{B(\frac{\tau + 2\rho}{3}) \setminus B(\rho)} |u|^2 dx + C (\tau - \rho)^{-1} \int_{B(\frac{\tau + 2\rho}{3}) \setminus B(\rho)} |u|^3 dx\\\nonumber
&&+C r (\tau - \rho)^{-1}  ||\nabla B||_{L^2(B(\frac{\tau + 2\rho}{3})\backslash B(\rho))} + C r^2 (\tau - \rho)^{-2} ||B - b||_{L^{6,\infty}(B(\frac{\tau + 2\rho}{3})\backslash B(\rho))}\nonumber\\
&&+C r^{\frac{4s-3}{s}}(\tau - \rho)^{-2} ||p - a||_{L^{\frac{3s}{3-s}}(B(\frac{\tau + 2\rho}{3})\backslash B(\rho))} +C r^{\frac{3s-3}{s}}(\tau - \rho)^{-1} ||\nabla p||_{L^s(B(\frac{\tau + 2\rho}{3})\backslash B(\rho))}\nonumber\\
\een

{\bf Step II. The bounded-ness of energy integral.}

Firstly, we estimate the first two terms of (\ref{ine:cacciopolli}) by following the estimate of \cite{Chae-Wolf2019}. Using $\nabla \cdot \mathbf{\Phi} = u$, integration by parts and H\"{o}lder inequality, we have
\beno
&& \int_{B(\frac{2\tau +\rho}{3})} |u|^3 \phi_2^3 dx \\
&=& \int_{B(\frac{2\tau +\rho}{3})} \partial_i (\Phi_{ij} - (\Phi_{ij})_{B(2r)}) u_j |u| \phi_2^3 dx \\
&=& - \int_{B(\frac{2\tau +\rho}{3})} (\Phi_{ij} - (\Phi_{ij})_{B(2r)}) \partial_i (u_j |u|) \phi_2^3 dx - \int_{B(\frac{2\tau +\rho}{3})} (\Phi_{ij} - (\Phi_{ij})_{B(2r)}) u_j |u| \partial_i (\phi_2^3) dx \\
&\leq& C\left(\int_{B(\frac{2\tau +\rho}{3})} |\nabla u|^2 dx\right)^{\frac 12} \left(\int_{B(\frac{2\tau +\rho}{3})} |\mathbf{\Phi} - (\mathbf{\Phi})_{B(2r)}|^6 dx\right)^{\frac 16} \left(\int_{B(\frac{2\tau +\rho}{3})} |u \phi_2|^3 dx\right)^{\frac 13}  \\
&+& C(\tau - \rho)^{-1} \left(\int_{B(\frac{2\tau +\rho}{3})} |\mathbf{\Phi} - (\mathbf{\Phi})_{B(2r)}|^3 dx\right)^{\frac 13} \left(\int_{B(\frac{2\tau +\rho}{3})} |u \phi_2|^3 dx\right)^{\frac 23},
\eeno
which implies from Young inequality that
\beno
\int_{B(\frac{2\tau +\rho}{3})} |u|^3 \phi_2^3 dx \leq C r^{\frac 34} \left(\int_{B(\frac{2\tau +\rho}{3})} |\nabla u|^2 dx\right)^{\frac 34} + C r^3 (\tau - \rho)^{-3}.
\eeno
Then we have
\ben\label{ine:I 2}\nonumber
(\tau - \rho)^{-1} ||u||_{L^3(B(\frac{\tau + 2\rho}{3}) \setminus B(\rho))}^3 &\leq& C r^{\frac 34} (\tau - \rho)^{-1} \left(\int_{B(\frac{2\tau +\rho}{3})} |\nabla u|^2 dx\right)^{\frac 34} \\
&& + C r^3 (\tau - \rho)^{-4}.
\een
By H\"{o}lder inequality, we get
\ben\label{ine:I 1}\nonumber
&&(\tau - \rho)^{-2} ||u||_{L^2(B(\frac{\tau + 2\rho}{3}) \setminus B(\rho))}^2 \leq C (\tau - \rho)^{-2} r ||u||_{L^3(B(\frac{\tau + 2\rho}{3}) \setminus B(\rho))}^2\\
&\leq& C r^{\frac 32} (\tau - \rho)^{-2} \left(\int_{B(\frac{2\tau +\rho}{3})} |\nabla u|^2 dx\right)^{\frac 12} + C r^3 (\tau - \rho)^{-4}.
\een


Let's substitute (\ref{ine:I 2}) and (\ref{ine:I 1}) into (\ref{ine:cacciopolli}). Note that $||B - b||_{L^{6,\infty}(B(2r))}\leq C$ from (\ref{eq:estimate b})  and by Young inequality we arrive at
\ben\label{ine:energy 1}
\int_{B(\rho)} |\nabla u|^2 + |\nabla B|^2 dx &\leq& \frac14 ||\nabla u||_{L^2(B(\frac{\tau + 2\rho}{3}))}^2+\frac14 ||\nabla B||_{L^2(B(\frac{\tau + 2\rho}{3}))}^2 \\ \nonumber
&& + C r^{\frac{3s-3}{s}}(\tau - \rho)^{-1} ||\nabla p||_{L^s(B(\frac{\tau + 2\rho}{3}))} + C r^4 (\tau - \rho)^{-4} \\ \nonumber
&& + C r^{\frac{4s-3}{s}}(\tau - \rho)^{-2} ||p - a||_{L^{\frac{3s}{3-s}}(B(\frac{\tau + 2\rho}{3}))}, \quad 1<s<3.
\een

Next we deal with the terms involving the pressure. By Lemma \ref{lem:Tsai's}, the estimate of (\ref{eq:stokes estimate}) ($n=3, q=s$ and $g=-u \cdot \nabla u + B \cdot \nabla B$) implies
\ben\label{eq:nabla p}
||\nabla p||_{L^s(B(\frac{\tau + 2\rho}{3}))} &\leq& C(n,s) \left(||u \cdot \nabla u - B \cdot \nabla B||_{L^s(B(\frac{2\tau +\rho}{3}))} +  (\tau-\rho)^{-2}||u||_{L^s(B(\frac{2\tau +\rho}{3}))}\right) \nonumber\\
&\doteq& C(n,s) \left(J_1 + J_2\right).
\een

On one hand, for $J_1$, using H\"{o}lder inequality we have
\beno
(J_1)^s &=& \int_{B(\frac{2\tau + \rho}{3})} |u \cdot \nabla u - B \cdot \nabla B|^s dx \\
&\leq& C \left(\int_{B(\frac{2\tau + \rho}{3})} |\nabla u|^2 dx\right)^{\frac s2} \left(\int_{B(\frac{2\tau + \rho}{3})} |u|^\frac{2s}{2-s} dx\right)^{\frac{2-s}{2}} \\
&&+ C \left(\int_{B(\frac{2\tau + \rho}{3})} |\nabla B|^2 dx\right)^{\frac s2} ||B||_{L^{6,\infty}(B(\frac{2\tau + \rho}{3}))}^s ||1||_{L^{\frac{3}{3-2s},\frac{2}{2-s}}(B(\frac{2\tau + \rho}{3}))} \\
&\leq& C \left(\int_{B(\frac{2\tau + \rho}{3})} |\nabla u|^2 dx\right)^{\frac s2} \left(\int_{B(\frac{2\tau + \rho}{3})} |u|^\frac{2s}{2-s} dx\right)^{\frac{2-s}{2}} \\
&&+ C r^{3-2s} \left(\int_{B(\frac{2\tau + \rho}{3})} |\nabla B|^2 dx\right)^{\frac s2}, \quad  1<s<\frac32.
\eeno
Due to the cut-off function (\ref{eq:psi 3}), we have
\beno
\int_{B(\frac{2\tau + \rho}{3})} |u|^\frac{2s}{2-s} dx \leq \int_{B(\tau)} |u \phi_3|^\frac{2s}{2-s} dx,
\eeno
hence,
\beno
&&\int_{B(\tau)} |u \phi_3|^\frac{2s}{2-s} dx \\
&=& \int_{B(\tau)} \partial_i(\Phi_{ij} - (\Phi_{ij})_{B(\tau)}) u_j |u|^\frac{4s - 4}{2-s} \phi_3^\frac{2s}{2-s} dx \\
&=& - \int_{B(\tau)}(\Phi_{ij} - (\Phi_{ij})_{B(\tau)}) \partial_i (u_j |u|^\frac{4s - 4}{2-s}) \phi_3^\frac{2s}{2-s} dx - \int_{B(\tau)}(\Phi_{ij} - (\Phi_{ij})_{B(\tau)}) u_j |u|^\frac{4s - 4}{2-s} \partial_i(\phi_3^\frac{2s}{2-s}) dx \\
&\leq& C \left(\int_{B(\tau)} |\nabla u|^2 dx\right)^{\frac 12} \left(\int_{B(\tau)} |u \phi_3|^\frac{2s}{2-s} dx\right)^{\frac{2s-2}{s}} \left(\int_{B(\tau)} |\mathbf{\Phi} - \mathbf{\Phi}_{B(\tau)}|^{\frac{2s}{4-3s}} dx\right)^{\frac{4-3s}{2s}} \\
&&+ (\tau - \rho)^{-1} \left(\int_{B(\tau)} |u \phi_3|^\frac{2s}{2-s} dx\right)^{\frac{3s-2}{2s}} \left(\int_{B(\tau)} |\mathbf{\Phi} - \mathbf{\Phi}_{B(\tau)}|^{\frac{2s}{2-s}} dx\right)^{\frac{2-s}{2s}}
\eeno
holds for all $s\in (1, \frac 43)$.
Using (\ref{eq: BMO property}) and Young inequality, we have
\ben\label{eq: u s}
\int_{B(\tau)} |u \phi_3|^\frac{2s}{2-s} dx \leq C r^{\frac{12-9s}{2(2-s)}} ||\nabla u||_{L^2(B(\tau))}^{\frac{s}{2-s}} + C r^3 (\tau - \rho)^{-\frac{2s}{2-s}}.
\een
Then we have
\beno
(J_1)^s &\leq& C r^{\frac{12-9s}{4}} ||\nabla u||_{L^2(B(\tau))}^{\frac{3s}{2}} + C r^{\frac{6-3s}{2}} (\tau - \rho)^{-s} ||\nabla u||_{L^2(B(\tau))}^s \\&&
+ C r^{3-2s} ||\nabla B||_{L^2(B(\tau))}^s.
\eeno
Using H\"{o}lder inequality and (\ref{eq: u s}), the term $J_2$ can be controlled by
\beno
(J_2)^s &=& (\tau - \rho)^{-2s} \int_{B(\frac{2\tau+\rho}{3})} |u|^s dx \\
&\leq& C r^{\frac{3s}{2}} (\tau - \rho)^{-2s} \left( \int_{B(\tau)} |u \phi_3|^\frac{2s}{2-s} dx \right)^{\frac {2-s}2} \\
&\leq& C r^{\frac{12-3s}{4}} (\tau - \rho)^{-2s} ||\nabla u||_{L^2(B(\tau))}^{\frac{s}{2}} + C r^3 (\tau - \rho)^{-3s},
\eeno
Combing the estimates of  $J_1$ and $J_2$, by (\ref{eq:nabla p}) we have
\ben\label{eq:stokes estimate'''} \nonumber
&&||\nabla p||_{L^s(B(\frac{\tau + 2\rho}{3}))} \\ \nonumber
&\leq& C r^{\frac{12-9s}{4s}} ||\nabla u||_{L^2(B(\tau))}^{\frac{3}{2}} + C r^{\frac{6-3s}{2s}} (\tau - \rho)^{-1} ||\nabla u||_{L^2(B(\tau))} \\ \nonumber
&&+ C r^{\frac{3-2s}{s}} ||\nabla B||_{L^2(B(\tau))} + C r^{\frac{12-3s}{4s}} (\tau - \rho)^{-2} ||\nabla u||_{L^2(B(\tau))}^{\frac{1}{2}} \\
&&+ C r^\frac{3}{s} (\tau - \rho)^{-3},  \quad  s\in (1, \frac 43).
\een
Furthermore, Poincar\'{e} inequality implies
\ben\label{eq:poincare}
||p - p_{B(\frac{\tau + 2\rho}{3})}||_{L^{\frac{3s}{3-s}}(B(\frac{\tau + 2\rho}{3}))} \leq C||\nabla p||_{L^s(B(\frac{\tau + 2\rho}{3}))}.
\een
Then choosing $a = p_{B(\frac{\tau + 2\rho}{3})}$ in  (\ref{ine:energy 1}), (\ref{eq:stokes estimate'''}) and (\ref{eq:poincare})  imply
\beno
&&\int_{B(\rho)} |\nabla u|^2 + |\nabla B|^2 dx \\&\leq& \frac 12 \int_{B(\tau)} |\nabla u|^2 + |\nabla B|^2 dx  + C r^4 (\tau - \rho)^{-4} \\
&& +C r^{\frac74}(\tau - \rho)^{-2}  ||\nabla u||_{L^2(B(\tau))}^{\frac{3}{2}} + Cr^{\frac52}(\tau - \rho)^{-3}  ||\nabla u||_{L^2(B(\tau))} \\ \nonumber
&&+ Cr^{2}(\tau - \rho)^{-2} ||\nabla B||_{L^2(B(\tau))} + Cr^{\frac{13}{4}}(\tau - \rho)^{-4}   ||\nabla u||_{L^2(B(\tau))}^{\frac{1}{2}}+ Cr^{4}(\tau - \rho)^{-5}
\eeno
Apply the iteration lemma (see, for example,  Lemma 3.1 in \cite{G83}), we have
\beno
\int_{B(\rho)} |\nabla u|^2 + |\nabla B|^2 dx &\leq& C r^{8} (\tau - \rho)^{-8}.
\eeno
Choosing $\rho = r$ and $\tau = 2r$, letting $r \rightarrow +\infty$ we have
\ben\label{ine:energy bounded}
\int_{\mathbb{R}^3} |\nabla u|^2 + |\nabla B|^2 dx \leq C.
\een

{\bf Step III. Vanishing property of Dirichlet integral.}

Next we prove the vanishing property of Dirichlet integral. It follows from (\ref{ine:energy bounded})  that $u,B \in \dot{W}^{1,2}(\mathbb{R}^3)$. Using Lemma \ref{lem:WZ} and the condition $u \in BMO^{-1}(\mathbb{R}^3)$, we have
\beno
|||u|^2||_{L^2(\mathbb{R}^3)} \leq C ||u||_{\dot{W}^{1,2}(\mathbb{R}^3)} ||u||_{BMO^{-1}(\mathbb{R}^3)},
\eeno
which means $u \in L^4(\mathbb{R}^3)$.
We claim that
\ben\label{ine:u L p}
u \in L^p(\mathbb{R}^3) \quad {\rm for ~~ all} \quad p \in [4,6].
\een
Indeed, for a new cut-off function $\zeta$ satisfying $\zeta = 1$ in $B(\frac 43 r)$ and $\zeta = 0$ on $B(2r)^c$, we get
\beno
\left(\int_{\mathbb{R}^3} |u \zeta|^6 dx\right)^\frac 13 \leq \int_{\mathbb{R}^3} |\nabla(u \zeta)|^2 dx \leq C\int_{B(2r)} |\nabla u|^2 dx + C r^{-\frac12} \left(\int_{B(2r)} |u|^4 dx\right)^\frac12,
\eeno
which implies $u \in L^6(\mathbb{R}^3)$. Similarly,
\ben\label{ine:b L p}
B \in L^6(\mathbb{R}^3),
\een
since
\beno
\left(\int_{\mathbb{R}^3} |B \zeta|^6 dx\right)^\frac 13  \leq C\int_{B(2r)} |\nabla B|^2 dx + C r^{-2} ||B||_{L^{6,\infty}(B(2r))}^2  ||1||_{L^{\frac32,1}(B(2r))}  <\infty.
\eeno

Recall the Cacciopolli inequality (\ref{ine:cacciopolli}), and
choose $\rho = r$, $\tau = 2r$, $s=\frac32$ and $b=B_{B(\frac{4}{3}r) \setminus B(r)}$. Then
\ben\label{ine:energy 2}
\int_{B(r)} |\nabla u|^2 + |\nabla B|^2 dx &\leq& C r^{-2} ||u||_{L^2(B(\frac43 r) \setminus B(r))}^2 + C r^{-1} ||u||_{L^3(B(\frac43 r) \setminus B(r))}^3  \nonumber \\
&&+ C ||B - b||_{L^{6,\infty}(B(\frac{4}{3}r)\backslash B(r))}  + C ||\nabla B||_{L^2(B(\frac{4}{3}r)\backslash B(r))} \nonumber\\
&& + C ||\nabla p||_{L^\frac 32 (B(\frac{4}{3}r)\backslash B(r))} + C ||p - a||_{L^{3}(B(\frac{4}{3}r)\backslash B(r))}.
\een
Using H\"{o}lder inequality, it follows  that
\beno
r^{-2} ||u||_{L^2(B(\frac43 r) \setminus B(r))}^2 \leq C r^{-\frac 12} ||u||_{L^4(B(\frac43 r) \setminus B(r))}^2\rightarrow 0,
\eeno
and
\beno
r^{-1} ||u||_{L^3(B(\frac43 r) \setminus B(r))}^3 \leq C r^{-\frac 14} ||u||_{L^4(B(\frac43 r) \setminus B(r))}^3\rightarrow 0,
\eeno
 as $r \rightarrow \infty$ since (\ref{ine:u L p}).
Moreover, note that
\beno
||B - b||_{L^6(B(\frac43 r) \setminus B(r))} \leq C ||\nabla B||_{L^2(B(\frac43 r) \setminus B(r))}\rightarrow 0  \quad {\rm as} \quad r \rightarrow \infty,
\eeno
which implies the third and the fourth terms are vanishing  as $r \rightarrow \infty$.

For the pressure, by (\ref{ine:b L p}) and (\ref{eq:stokes estimate}) in Lemma \ref{lem:Tsai's},  we have
\beno
||\nabla p||_{L^\frac32(B(\frac 43))}
&\leq &C   ||u||_{L^6(B(2r))} ||\nabla u||_{L^2(B(2r))} +  C||B||_{L^6(B(2r))} ||\nabla B||_{L^2(B(2r))} \nonumber\\
&& +Cr^{-1/2}  ||u||_{L^6(B(2r))}.
\eeno
Letting $r \rightarrow \infty$, we obtain $\nabla p \in L^\frac 32(\mathbb{R}^3)$, which yields $||\nabla p||_{L^\frac 32 (B(\frac{4}{3}r)\backslash B(r))} \rightarrow 0$.
Choosing $a = p_{B(\frac{4}{3}r)\backslash B(r)}$, it immediately implies
\beno
||p - a||_{L^{3}(B(\frac{4}{3}r)\backslash B(r))} \leq ||\nabla p||_{L^\frac 32 (B(\frac{4}{3}r)\backslash B(r))} \rightarrow 0.
\eeno

Hence, (\ref{ine:energy 2}) implies
\beno
\int_{B(r)} |\nabla u|^2 + |\nabla B|^2 dx \rightarrow 0 \quad {\rm as} \quad r \rightarrow \infty,
\eeno
which means $u \equiv 0$ and $B \equiv 0$, since $u \in L^{6}(\mathbb{R}^3)$ and $B \in L^{6,\infty}(\mathbb{R}^3)$. The proof is complete.
\end{proof}

\section{Proof of Theorem \ref{thm:viscidity}}

\begin{proof}
 Assume that $\eta(x) \in C_0^{\infty}(\mathbb{R}^3)$ with $ 0 \leq \eta \leq 1$, satisfying that\\
i) for $r > 1$, $r \leq \rho < \tau \leq 2r$,
\begin{equation*} \eta(x)=\eta(|x|)=
\left\{\begin{array}{llll}
1,\quad ~ {\rm in} ~ B(\rho);\\
0, \quad ~ {\rm in} ~ \mathbb{R}^3 \backslash B(\tau)
\end{array}\right.
\end{equation*}
ii)
\beno
&&|\nabla \eta| \leq \frac{C}{\tau - \rho},\quad |\nabla^2 \eta| \leq \frac{C}{(\tau - \rho)^2}.
\eeno
Since $B(\tau)$ is a star-like domain(see P38, \cite{Galdi}), due to Theorem III 3.1 in \cite{Galdi}, there exists a constant $C(s)$ and a vector-valued function $w : B(\tau) \rightarrow \mathbb{R}^3$ such that
$w \in W^{1,s}_0(B(\tau))$, and
$\nabla \cdot w = u \cdot \nabla \eta$
satisfying
\ben\label{esti-ws'}
\int_{B(\tau)} |\nabla w(x)|^s \,dx dy
\leq C(s) \int_{B(\tau)} |u \cdot \nabla \eta|^s \,dx, \quad 1 < s < \infty.
\een
Then  multiplying $(\ref{eq:VMHD})_1$ by $u \eta-w $, integrating it over $B(2r)$ and integration by parts yield that
\ben\label{eq:energy u} \nonumber
&&\lambda_1 \int_{B(2r)} |\partial_1 u|^2 \eta dx + \lambda_2 \int_{B(2r)} |\partial_2 u|^2 \eta dx + \lambda_3 \int_{B(2r)} |\partial_3 u|^2 \eta dx \\ \nonumber
&=& \lambda_1 \int_{B(2r)} \partial_1 u \cdot \partial_1 w dx + \lambda_2 \int_{B(2r)} \partial_2 u \cdot \partial_2 w dx + \lambda_3 \int_{B(2r)} \partial_3 u \cdot \partial_3 w dx \\ \nonumber
&+& \frac 12 \lambda_1 \int_{B(2r)} |u|^2 \partial_{11} \eta dx + \frac 12 \lambda_2 \int_{B(2r)} |u|^2 \partial_{22} \eta dx + \frac 12 \lambda_3 \int_{B(2r)} |u|^2 \partial_{33} \eta dx \\ \nonumber
&-& \int_{B(2r)} u \cdot \nabla w \cdot u dx + \frac 12 \int_{B(2r)} |u|^2 u \cdot \nabla \eta dx + \int_{B(2r)} B \cdot \nabla w \cdot B dx \\
&+& \int_{B(2r)} B \cdot \nabla B \cdot u \eta dx.
\een
Similarly, by multiplying $(\ref{eq:VMHD})_2$ by $B \eta$ we derive
\ben\label{eq:energy B}
&&\mu_1 \int_{B(2r)} |\partial_1 B|^2 \eta dx + \mu_2 \int_{B(2r)} |\partial_2 B|^2 \eta dx + \mu_3 \int_{B(2r)} |\partial_3 B|^2 \eta dx \\ \nonumber
&=& \frac 12 \mu_1 \int_{B(2r)} |B|^2 \partial_{11} \eta dx + \frac 12 \mu_2 \int_{B(2r)} |B|^2 \partial_{22} \eta dx + \frac 12 \mu_3 \int_{B(2r)} |B|^2 \partial_{33} \eta dx \\ \nonumber
&+& \frac 12 \int_{B(2r)} |B|^2 u \cdot \nabla \eta dx - \int_{B(2r)} B \cdot \nabla B \cdot u \eta dx - \int_{B(2r)} B \cdot u (B \cdot \nabla) \eta dx.
\een
Letting
\beno
G(r) = \sum_{i=1}^3 \lambda_i \int_{B(r)} |\partial_i u|^2 dx + \sum_{i=1}^3 \mu_i \int_{B(r)} |\partial_i B|^2 dx,
\eeno
and taking the sum of  (\ref{eq:energy u}) and (\ref{eq:energy B}), we have
\beno
G(\rho) &\leq& \lambda_1 \int_{B(2r)} \partial_1 u \cdot \partial_1 w dx + \lambda_2 \int_{B(2r)} \partial_2 u \cdot \partial_2 w dx + \lambda_3 \int_{B(2r)} \partial_3 u \cdot \partial_3 w dx \\ \nonumber
&+& \frac 12 \lambda_1 \int_{B(2r)} |u|^2 \partial_{11} \eta dx + \frac 12 \lambda_2 \int_{B(2r)} |u|^2 \partial_{22} \eta dx + \frac 12 \lambda_3 \int_{B(2r)} |u|^2 \partial_{33} \eta dx \\ \nonumber
&-& \int_{B(2r)} u \cdot \nabla w \cdot u dx + \frac 12 \int_{B(2r)} |u|^2 u \cdot \nabla \eta dx + \int_{B(2r)} B \cdot \nabla w \cdot B dx \\
&+& \frac 12 \mu_1 \int_{B(2r)} |B|^2 \partial_{11} \eta dx + \frac 12 \mu_2 \int_{B(2r)} |B|^2 \partial_{22} \eta dx + \frac 12 \mu_3 \int_{B(2r)} |B|^2 \partial_{33} \eta dx \\ \nonumber
&+& \frac 12 \int_{B(2r)} |B|^2 u \cdot \nabla \eta dx - \int_{B(2r)} B \cdot u (B \cdot \nabla) \eta dx \doteq\sum_{i=1}^{14} I_i.
\eeno

Since the terms $I_1$, $I_2$ and $I_3$ are similar, we compute the first term. Letting $T(\rho,\tau) = B(\tau) \setminus B(\rho)$, using H\"{o}lder and Young inequality, (\ref{esti-ws'}) implies
\beno
|I_1| &\leq& \lambda_1 \int_{B(\tau)} |\partial_1 u| |\partial_1 w| dx \\
&\leq& \varepsilon \lambda_1 \int_{B(\tau)} |\partial_1 u|^2 dx + C \lambda_1 \int_{B(\tau)} |\partial_1 w|^2 dx \\
&\leq& \varepsilon G(\tau) + C  \lambda_1 (\tau - \rho)^{-2} \int_{T(\rho,\tau)} |u|^2 dx \\
&\leq& \varepsilon G(\tau) + C  \lambda_1 r (\tau - \rho)^{-2} ||u||_{L^3(T(\rho,\tau))}^2,
\eeno
where $\varepsilon>0$, to be decided. Similarly,
\beno
|I_2|+|I_3| \leq \varepsilon G(\tau) + C(\lambda_2+ \lambda_3) r (\tau - \rho)^{-2} ||u||_{L^3(T(\rho,\tau))}^2.
\eeno

The estimates of $I_4$, $I_5$ and $I_6$ can be obtained directly from H\"{o}lder's inequality:
\beno
|I_4|+ |I_5|+|I_6|\leq  C (\lambda_1+ \lambda_2+\lambda_3) r (\tau - \rho)^{-2} ||u||_{L^3(T(\rho,\tau))}^2.
\eeno
Similarly, the terms of $I_{10}$, $I_{11}$ and $I_{12}$ are also obtained:
\beno
|I_{10}|+ |I_{11}|+|I_{12}|\leq C ( \mu_1+\mu_2+\mu_3) r^2 (\tau - \rho)^{-2} ||B||_{L^6(T(\rho,\tau))}^2.
\eeno

For the term $I_7$, using (\ref{esti-ws'}) again, we have
\beno
|I_7| &\leq& \int_{B(2r)} |u|^2 |\nabla w| dx \leq \left(\int_{B(\tau)} |u|^3 dx\right)^{\frac 23} \left(\int_{B(\tau)} |\nabla w|^3 dx\right)^{\frac 13} \\
&\leq& C (\tau - \rho)^{-1} ||u||_{L^3(T(\rho,\tau))}||u||_{L^3(B(\tau))}^2,
\eeno
and
\beno
|I_8| \leq C (\tau - \rho)^{-1} ||u||_{L^3(T(\rho,\tau))}||u||_{L^3(B(\tau))}^2.
\eeno
Similarly,
\beno
|I_9| &\leq& \int_{B(2r)} |B|^2 |\nabla w| dx \leq C r ||B||_{L^6(B(\tau))}^2 \left(\int_{B(\tau)} |\nabla w|^3 dx\right)^{\frac 13} \\
&\leq& C r (\tau - \rho)^{-1} ||B||_{L^6(B(\tau))}^2 ||u||_{L^3(T(\rho,\tau))}.
\eeno
The terms $I_{13}$ and $I_{14}$ can be treated as the term $I_9$:
\beno
|I_{13}| + |I_{14}|  \leq C r (\tau - \rho)^{-1} ||B||_{L^6(T(\rho,\tau))}^2 ||u||_{L^3(T(\rho,\tau))}.
\eeno

Collecting all the terms of $I_i$, we arrive at
\beno
G(\rho) &\leq& 2\varepsilon G(\tau) +  C  (\lambda_1+\lambda_2 +\lambda_3) r (\tau - \rho)^{-2} ||u||_{L^3(T(\rho,\tau))}^2\\
 &&+ C (\mu_1+\mu_2+\mu_3) r^2 (\tau - \rho)^{-2} ||B||_{L^6(T(\rho,\tau))}^2 \\
&&+ C (\tau - \rho)^{-1} ||u||_{L^3(T(\rho,\tau))}||u||_{L^3(B(\tau))}^2 + C r (\tau - \rho)^{-1} ||B||_{L^6(B(\tau))}^2 ||u||_{L^3(T(\rho,\tau))}\\
&\leq& 2\varepsilon G(\tau) +  C  (\lambda_1+\lambda_2 +\lambda_3) r (\tau - \rho)^{-2} ||u||_{L^3(T(r,2r))}^2\\
 &&+ C (\mu_1+\mu_2+\mu_3) r^2 (\tau - \rho)^{-2} ||B||_{L^6(T(r,2r))}^2 \\
&&+ C (\tau - \rho)^{-1} ||u||_{L^3(T(r,2r))}||u||_{L^3(B(\tau))}^2 + C r (\tau - \rho)^{-1} ||B||_{L^6(B(\tau))}^2 ||u||_{L^3(T(r,2r))}.
\eeno
Choosing $\varepsilon = \frac 14$ and using the iteration in \cite{G83}, we have
\beno
G(\rho) &\leq& C  (\lambda_1+\lambda_2 +\lambda_3) r (\tau - \rho)^{-2} ||u||_{L^3(T(r,2r))}^2\\
 &&+ C (\mu_1+\mu_2+\mu_3) r^2 (\tau - \rho)^{-2} ||B||_{L^6(T(r,2r))}^2 \\
&&+ C (\tau - \rho)^{-1} ||u||_{L^3(T(r,2r))}||u||_{L^3(B(\tau))}^2 + C r (\tau - \rho)^{-1} ||B||_{L^6(B(\tau))}^2 ||u||_{L^3(T(r,2r))}.
\eeno
Taking $\rho = r$ and $\tau = 2r$,   letting $r \rightarrow \infty$, we have
\beno
\sum_{i=1}^3 \lambda_i \int_{\mathbb{R}^3} |\partial_i u|^2 dx + \sum_{i=1}^3 \mu_i \int_{\mathbb{R}^3} |\partial_i B|^2 dx = 0,
\eeno
since $u \in L^3(\mathbb{R}^3)$ and $B \in L^6(\mathbb{R}^3)$. Then there exist $i \in \{1,2,3\}$ and $j \in \{1,2,3\}$, such that $\partial_i u, \partial_j B \equiv 0$ due to the known condition (\ref{eq:condition}). Without loss of generality, we assume that $\partial_1 u = \partial_1 B = 0$, which means $u(x_1, x_2, x_3) = u(x_2, x_3)$ and $B(x_1, x_2, x_3) = B(x_2, x_3)$. Then
\beno
\int_{\mathbb{R}^2} \int_{-\infty}^{+\infty} |u(x_2, x_3)|^3dx_1 dx_2 dx_3 =\infty,
\eeno
if $u$ is a nontrivial solution, which is a contradiction with $u \in L^3(\mathbb{R}^3)$. Thus we have $u\equiv 0$ in $\mathbb{R}^3$. Using the same idea, we have $B \equiv 0$ in $\mathbb{R}^3$. The proof is complete.
\end{proof}

\noindent {\bf Acknowledgments.}
W. Wang was supported by NSFC under grant 12071054 and 11671067.


\end{document}